\documentclass{elsart}

\usepackage{latexsym, epsfig, rotfloat}
\begin{document}

\newcommand {\heading} [1]
{\par \vspace{5mm} \noindent
{\Large \bf #1 } \par \nopagebreak \vspace{5mm} \nopagebreak}

%\hyphenpenalty=10000
%\sloppy

\newcommand {\expand} {{\sc expand}}
\newcommand {\Aaa} {{\mbox{\scriptsize $A_{11}$}}}
\newcommand {\Aab} {{\mbox{\scriptsize $A_{12}$}}}
\newcommand {\Aba} {{\mbox{\scriptsize $A_{21}$}}}
\newcommand {\lra} {\Leftrightarrow}
\newcommand {\beq} {\vspace{-4mm}\begin{eqnarray}}
\newcommand {\eeq} {\end{eqnarray}}
\newcommand {\beqa} {\vspace{-4mm}\begin{eqnarray}}
\newcommand {\eeqa} {\end{eqnarray}}
\newcommand {\beqas} {\vspace{-4mm}\begin{eqnarray*}}
\newcommand {\eeqas} {\end{eqnarray*}}
\newcommand {\beqs} {\vspace{-4mm}\begin{eqnarray*}}
\newcommand {\eeqs} {\end{eqnarray*}}
\newcommand {\bpd} {\left\{ \begin{array}{ll}}
\newcommand {\epd} {\end{array} \right.}
\newcommand{\nn}{\nonumber}
\newcommand{\dis}{\displaystyle}
\def\QED{\hfill{\bf Q.E.D.}}
%\newcommand {\lb}[1] {\mbox{\b{$#1$}}}
  % To use: \lb{X} will make an underlined MATHEMATICAL X when used either
  % within math mode or within text.
%\newcommand {\ub}[1] {\mbox{$\bar{#1}$}}
%BLACKBOARD BOLD CHARACTERS        R, Z, C, Q, N:
\newcommand{\R}{\mbox{$I\!\!R$}}
\newcommand{\Z}{\mbox{$Z\!\!\!Z$}}
\newcommand{\C}
        {\mbox{$C\!\!\!{\vrule height1.5ex width1pt depth-0.12ex}\,\,\,$}}
\newcommand{\Q}
        {\mbox{$Q\!\!\!{\vrule height1.5ex width1pt depth0pt}\,\,\,$}}
\newcommand{\N}{\mbox{$I\!\!N$}}          % N
\newcommand {\cg}[2]
{{#1}{}}
%{\marginpar{x}{#1}{}}
% marginal note/proposed text/old text

%\newtheorem{thm}{Theorem}
%\newtheorem{lem}{Lemma}
%\newtheorem{prop}{Proposition}
%\setcounter{page}{0}

%\input{cover}

\begin{frontmatter}
\title{The simplest examples where the simplex method cycles \\
and conditions where \expand~fails to prevent cycling\thanksref{EPSRC}}
\author{J.A.J. Hall and %\thanksref{jajh},
%},
%\author{
K.I.M. {McKinnon}\thanksref{both}}
\address{Dept.\ of Mathematics and Statistics, 
               University of Edinburgh, EH9 3JZ,  UK}
\thanks[EPSRC]{Work supported by EPSRC grant GR/J0842}
%\thanks[jajh]{jajhall@maths.ed.ac.uk}
\thanks[both]{jajhall@maths.ed.ac.uk,ken@maths.ed.ac.uk}
\begin{abstract}
This paper introduces a class of linear programming examples which
cause the simplex method to cycle indefinitely and which are the
simplest possible examples showing this behaviour. The structure of
examples from this class repeats after two iterations. Cycling is shown
to occur for both the most negative reduced cost and steepest edge
column selection criteria. In addition it is shown that the
\expand~anti-cycling procedure of Gill
\emph{et al.}\ is not guaranteed to prevent cycling.

\end{abstract}
\end{frontmatter}

\section{Introduction}

Degeneracy in linear programming is of both theoretical and practical
importance. It occurs whenever one or more of the basic variables is
at its bound. An iteration of
the simplex method may then fail to improve the objective function.  The
simple proof of finiteness of the simplex algorithm relies on a strict
improvement in the objective function at each iteration and the fact
that the simplex method visits only basic solutions, of which there
is a finite number. However if the problem is degenerate
there is the possibility of a consecutive sequence of
iterations occurring with no change in the objective function and with
the eventual return to a previously encountered basis. Examples such
as Beale's \cite{beal55} have been constructed to
show that this can happen, though such examples do seem to be very
rare in practice. A more common practical situation is where a long
but finite sequence of iterations occurs without the objective
function improving---a situation called {\em stalling}---and this
can degrade the algorithm's performance.

A related issue is the behaviour of the simplex algorithm in
the presence of roundoff error. At a degenerate vertex there is
a serious danger of selecting pivots that are small and have a
high relative error.

A wide range of methods have been suggested to avoid these problems.
\begin{description}
\item[
Lexicographic ordering:] These methods are guaranteed to terminate
in exact
arithmetic but are often prohibitively expensive to implement for the
revised simplex method and do not address the problem of inexact
arithmetic.
\item[Primal-dual alternation:] These methods were introduced by
Balinski and Gomory \cite{BalGom63} and have recently been developed by
Fletcher \cite{Fl88,FlHa90,FL93}. Some of these methods
guarantee to terminate in exact arithmetic and also exhibit good behaviour
with inexact arithmetic.
\item[ Constraint perturbation and feasible set
enlargement:] These methods attempt to reduce the likelyhood of
cycling and also attempt to improve the numerical behaviour and reduce
the number of iterations. The Devex and \expand~procedures described
below are of this type. In addition it is claimed that
stalling cannot occur with \expand\ \cg{with exact arithmetic}{}.
Wolfe's method is a recursive perturbation method which
guarantees termination in exact arithmetic.
\end{description}
%\end{itemize}

In \cite{Wf63} Wolfe introduced a perturbation method which is
guaranteed to terminate in a finite number of steps in exact arithmetic.
In this method, whenever a degenerate vertex is encountered, the bounds
producing the degeneracy are expanded in such a way that the current
vertex is no longer degenerate.  Other bounds on the basic variables
are temporarily ignored. The simplex method works on this modified
problem until an unbounded direction is found.  If the bound expansion
is random, it is highly unlikely that further degenerate vertices will
be encountered before the unbounded direction is found. However if a
further degenerate vertex is discovered, it is guaranteed to have
fewer active constraints. The perturbation process is repeated and
\cg{}{since it is guaranteed that each successive degenerate vertex has
fewer active constraints,}
after a finite number of steps a
non-degenerate vertex is reached with an unbounded direction. This
direction is then used in the original problem to give an edge leading
out of the degenerate vertex. It is not obvious how to extend this
method to the case of inexact arithmetic, as there is then no
obvious criterion for what constitutes a degenerate vertex.

In \cite{Harris} Harris introduced the Devex row selection method,
which allowed small violations of the constraints and used the
resulting flexibility to choose the largest pivot. This has the
advantage of both avoiding unnecessarily small pivots and reducing the
number of iterations. The disadvantage is that the constraints are
violated and some steps are \cg{negative}{in the negative direction}.
The variable leaving the basis does not normally do so at one of its
bounds, but is shifted to that value, resulting in inconsistent values
for the basic variables.  The method attempts to correct this
inconsistency at regular intervals (usually after each reinversion) by
doing a {\em reset}, in which the basic variable values are
recalculated from the values of the nonbasic variables. This can
produce infeasible values for the basic variables (i.e. outside the
specified tolerance) so there is no guarantee that progress has been
made. However the method seems to be effective in practice in reducing
the number of iterations taken, and variants of it are used in some
commercial codes.

Gill {\em et al} \cite{GMSW89} developed the \expand~method in an
attempt to improve on the good features of the Devex method of Harris
and also to incorporate some features of Wolfe's method which
guarantee finite termination. The performance of {\sc minos} was
significantly improved by the incorporation of \expand. At each
iteration of the
\expand~method the
bounds are expanded by a small amount. As in Devex, the largest pivot
that does not lead to any constraint violation (beyond the current
expanded position) is chosen. If the normal step for the largest pivot
is sufficiently positive, it is taken; otherwise a small positive step
is taken. In all cases the variable values stay within their expanded
bounds. Because at every iteration the nonbasic variable is moved a
positive amount in the direction that improves the objective
function, \cg{the objective}{it} can never return to a previous value\cg{
so no previous solution can recur}{}.
% \cg{Leave to later?}{}{However, since the
%bounds are being expanded at every iteration, the optimal solution to
%the expanded problem may improve faster than the current solution, so
%there is no guarantee of reaching the optimal solution.}
 
In this paper we introduce and analyse the simplest possible
class of cycling examples, the {2/6-cycle} class.
In Section \ref{IntroEx} we present an example of
this class which cycles when using the most negative reduced cost column
selection criterion.  In Section \ref{2/6form} the general form of
such examples is derived.  In Section \ref{SteepestEdge} a variation
of the example is introduced which cycles for the steepest-edge column
selection rule. In Section
\ref{EXPAND} the behaviour of the \expand~procedure is analysed and a
simple necessary and sufficient condition is derived for indefinite
cycling to occur\cg{}{when using expand}.

\section{Introductory example} \label{IntroEx}

We first solve the four variable, two constraint problem
(\ref{prob1}) by the simplex method. The analysis later in the paper
shows how to derive examples of this form. The problem is unbounded. A
bounded example with identical behaviour can be obtained by adding the
upper bound constraints $x_1\leq 1$ and $x_2\leq 1$, either as implicit
upper bounds or with one or more explicit constraints. The variable
to enter the basis will be chosen by the most negative reduced cost
criterion and, where there is a tie for the variable to leave the
basis, the variable in the row with the largest pivot will be chosen.
\beqa
%\mbox{Max} \;\;\;\; I = x_1 - 5.5 x_2 + 0.75 x_3 - 5.75 x_4 \nonumber \\
%\mbox{Subject to} \;\;\;\;
%           2.5 x_1 - 19.5 x_2 -3.5 x_3 + 19.5 x_4 & \leq & 0 \label{prob1} \\
%           0.5 x_1 -  3.5 x_2 -0.5 x_3 + 2.5 x_4 & \leq & 0 \nonumber \\
%                 x_i & \geq & 0, \;\; i = 1 \ldots 4 \nonumber
%\eeqa
\mbox{Max} \;\;\;\; z = 2.3 x_1 + 2.15 x_2 -13.55 x_3 -0.4 x_4,& & \nonumber \\
\mbox{subject to} \;\;\;\;\;\;\;\;
           0.4 x_1 +  0.2 x_2 -1.4 x_3 -  0.2 x_4 & \leq & 0, \label{prob1} \\
          -7.8 x_1 -  1.4 x_2 +7.8 x_3 + 0.4 x_4 & \leq & 0, \nonumber \\
                 x_j & \geq & 0, \;\; j = 1 \ldots 4. \nonumber
\eeqa
After introducing slack variables $x_5$ and $x_6$ and writing the
equations in detached coefficient form we get tableau $T^{(1)}$.
All the variables are initially zero and will remain zero at every iteration.
In the first iteration $x_1$ is chosen to 
enter the basis. There is only one positive entry in the $x_1$ column, so
there is a unique pivot choice with $x_5$ leaving the basis. This
basis change leads to tableau $T^{(2)}$.
In the second iteration $x_2$ is
chosen to enter the basis. In the normal ratio test there is a tie
between $x_6$ and $x_1$ to leave the basis. Breaking the tie by
using the larger pivot
(as is normal for numerical stability) gives $x_6$ to leave the basis, and
the basis change yields tableau $T^{(3)}$.

\vspace{4mm}
\centerline{
\begin{tabular}{rrrrrrrcrr}
$x_1$ & $x_2$ & $x_3$ & $x_4$ & $x_5$ & $x_6$ & $z$ \\
\hline
   0.4 &   0.2 &  -1.4 &  -0.2 &  1.0  &     &     & = & 0 \\
  -7.8 &  -1.4 &   7.8 &   0.4 &       & 1.0 &     & = & 0 & ~~~~$T^{(1)}$ \\
  -2.3 &  -2.15 & 13.55 &  0.4 &       &     & 1.0 & = & 0 \\
\hline
   1.0 &   0.5  &  -3.5 &  -0.5 &   2.5 &     &     &= & 0 \\
       &   2.5  & -19.5 &  -3.5 &  19.5 & 1.0 &  &= & 0 & ~~~~~~~~$T^{(2)}$ \\
       &  -1.0  &   5.5 & -0.75 &  5.75 &     & 1.0 & = & 0 \\
\hline
   1.0 &     &  0.4 &   0.2 &  -1.4 &  -0.2 &     & = & 0 \\
       & 1.0 & -7.8 &  -1.4 &   7.8 &   0.4 &     & = & 0 & ~~~~$T^{(3)}$ \\
       &     & -2.3 &  -2.15 & 13.55 &  0.4 & 1.0 & = & 0 \\
\hline
\end{tabular}
}
%\vspace{4mm}

Note that tableau $T^{(3)}$ is the same as tableau $T^{(1)}$ with the $x$
variable columns shifted cyclically two columns to the right. It follows
that this example will return to tableau $T^{(1)}$ after a further 4
iterations and therefore will cycle indefinitely\cg{ with a cycle length of 6}{}. In this example there
are only two sets of coefficients: $T^{(3)}$ and $T^{(5)}$ are the same as $T^{(1)}$
with the $x$ variable columns shifted cyclically 2 and 4 columns to the
right, and $T^{(4)}$ and $T^{(6)}$ are the same as $T^{(2)}$ again shifted
cyclically 2 and 4 columns to the right. We refer to such examples
as 2/6-cycle examples. In this paper we restrict attention to
2/6-cycle examples as they are more elegant and easier to analyse than
6/6-cycle examples, such as Beale's example, which take 6 iterations to
repeat  the same coefficients. \cg{}{
??It is possible to obtain
6/6-cycle examples with similar behaviour 
by perturbing the 2/6-cycle examples.??}
All the results here are demonstrated for
2/6-cycle examples. However the 2/6 property is not needed for the results
and indeed 6/6-cycle examples can be formed by perturbing the 2/6-cycle
examples given in this paper.

\section{The form of 2/6-cycle examples} \label{2/6form}

The following analysis was used to construct the above example.
Let the $3 \times 6$ matrix $M^{(1)}$ be formed from the $x$ columns
of $T^{(1)}$ as follows
\beqs
M^{(1)} =
\left[ \begin{tabular}{ccc}
$A$ & $B$ & $I$ \\
$a$ & $b$ & $0$ 
\end{tabular}
\right],
\eeqs
where $A$, $B$ and $I$ are $2 \times 2$ blocks of the constraint rows
and $a$, $b$ and $0$ are $1 \times 2$ blocks of the objective row. To
be able to pivot on the (1,1) and (2,2) entries in iterations 1 and 2,
we require $A$ to be non-singular. These pivoting operations
yield
tableau $T^{(3)}$, whose submatrix formed from the $x$ columns has the form
\beqs
M^{(3)}=
\left[ \begin{tabular}{ccc}
$I$ & $A^{-1}B$ & $A^{-1}$ \\
$0$ & $b-aA^{-1}B$ & $-aA^{-1}$
\end{tabular}
\right].
\eeqs
For the constraint pattern to repeat after these two iterations we
require $A = A^{-1}B$ and $B = A^{-1}$, which occurs if and only if
$A^3 = I$.
This implies that the eigenvalues, $\lambda$, of $A$ satisfy
\beq
 \lambda^3  =  1 
 \iff (\lambda^2 + \lambda+1)(\lambda-1) & = & 0. \label{fulllambda}
\eeq
For a $2 \times 2$ real matrix $A$ there must either be 2 real
eigenvalues or a complex conjugate pair.

It follows from (\ref{fulllambda}) that if $A$ has real eigenvalues
they must both have the value 1, in which case the $2\times2$ matrix
polynomial $A^2+A+I$ has two real eigenvalues of 3 and is therefore
non-singular. Since $(A-I)(A^2+A+I) = A^3-I=0$, it follows that $A=I$
in this case. It is then easy to show that $a=b=0$, which is of no
interest as it corresponds to a zero cost row.

The other possibility is that $A$ has
a complex conjugate pair of eigenvalues, and it follows from (\ref{fulllambda})
that they must satisfy
\beq
   \lambda^2 + \lambda+1 = 0. \label{chareqn}
\eeq
The characteristic equation of a general $2 \times 2$ matrix $A$ is
\beq
   \lambda^2 - (A_{11}+A_{22})\lambda+(A_{11}A_{22}-A_{21}A_{12}) = 0.
   \label{chareqn2}
\eeq
Equations (\ref{chareqn}) and (\ref{chareqn2})
hold for the two distinct
values of $\lambda$,
so \cg{}{it follows that} for a suitable 2/6-cycle example we require
$A_{11}+A_{22}  =  -1$ and
$A_{11}A_{22}-A_{21}A_{12}  =  1$.
From these it follows that
\beq
   -A_{21}A_{12} = 1+A_{11} +A^2_{11}. \label{Aelrelation}
\eeq
Conversely, any $2\times 2$ matrix such that $A_{11}+A_{22}=-1$
and~(\ref{Aelrelation}) holds has characteristic equation~(\ref{chareqn}).
Since a matrix satisfies its own characteristic equation,
$A^2+A+I=0$, from which it follows that $A^3=I$.

The objective function will repeat after 2 iterations if and only if
$b-aA^{-1}B = a$ and $b=-aA^{-1}$. This occurs if and only if
$a(A^2+A+I) = 0$, which holds for all $a$
%as $A$ must satisfy its own characteristic equation
since $A^2+A+I=0$. There is therefore no restriction
on $a$. Since the scaling of the objective row is
arbitrary we take $a$ to have the form
\beqs
a = [-1,\mu],
\eeqs
where there is no restriction on the value of $\mu$. It follows that
there is a three parameter family of 2/6-cycle examples: the parameters
can be chosen as $\mu$, $A_{11}$ and $A_{12}$.

For arbitrary $a$, the vector $b$ must satisfy
\beq
   b = -aA^{-1}. \label{baArel}
\eeq
Since $A$ is real and $A^3=I$, $\det(A)=1$. Hence
\beq
   B = A^{-1} = \left[
\begin{tabular}{cc}
    $-(A_{11}+1)$ & $-A_{12}$ \\
    $-A_{21}$    & $A_{11}$
\end{tabular}
\right], \nonumber
\eeq
\beq
   \mbox{and} \;\;\; b = [-(A_{11}+1) + \mu A_{21}, -A_{12}-\mu A_{11}],
   \nonumber
\eeq
and
follows that the general form of $M^{(1)}$ and $M^{(2)}$ for the 2/6-cycle
examples with the pivot sequence fixed is as in Table \ref{tab:CoefVals}.

\begin{table}
\caption{Coefficient values over two iterations for 2/6-cycle examples}
%\vspace{4mm}
%\hspace{-3mm}
{\scriptsize
%{\footnotesize
\centerline{\begin{tabular}{ccccccc}
&$x_1$ & $x_2$ & $x_3$ & $x_4$ & $x_5$ & $x_6$ \\ \noalign{\smallskip}
\hline
& $\Aaa$  & $\Aab$ & $-(\Aaa+1)$ & $-\Aab$ &   1  \\  \noalign{\smallskip}
 $M^{(1)}=$& $\Aba$  & $-(\Aaa+1)$ &  $-\Aba$ & $\Aaa$ &       &   1 \\  \noalign{\smallskip}
& -1  & $\mu$ & $-(\Aaa+1)+\mu \Aba$ & $-\Aab-\mu \Aaa$ \\  \noalign{\smallskip}
\hline
& 1 & $\frac{A_{12}}{A_{11}}$ & $-(1+\frac{1}{A_{11}})$ &
$-\frac{A_{12}}{A_{11}}$ & $\frac{1}{A_{11}}$ \\  \noalign{\smallskip}
$M^{(2)}=$ & & $\frac{1}{A_{11}}$ & $\frac{\Aba}{A_{11}}$ & $-(1+\frac{1}{A_{11}})$ &
 $-\frac{A_{21}}{A_{11}}$ & 1 \\  \noalign{\smallskip}
& & $\mu+\frac{A_{12}}{A_{11}}$ &
  $\mu \Aba-(2+\Aaa+\frac{1}{A_{11}})$ &
  $-\Aab(1+\frac{1}{A_{11}})-\mu \Aaa$ & $\frac{1}{A_{11}}$
\vspace{2mm}
 \\
\hline
\end{tabular}
}
}
\vspace{4mm}
\label{tab:CoefVals}
\end{table}

Proposition \ref{Prop26pattern} summarises these results.

\begin{prop} \label{Prop26pattern}
Assume the cost row is nonzero and the 2/6-cycle pattern of pivots is
selected. Then the necessary and sufficient conditions for the
coefficient pattern to repeat after two iterations are that the
coefficients have the form given in tableau $M^{(1)}$ of Table
\ref{tab:CoefVals}, and that $A_{11}$, $A_{21}$ and $A_{12}$ satisfy
(\ref{Aelrelation}).
\end{prop}

We now deduce the inequality relations that must be satisfied for
the simplex method to select (1,1) and (2,2) as pivot elements.
In order for (1,1) to be a pivot in tableau $M^{(1)}$ we require
\beq
   A_{11} > 0. \label{A11pos}
\eeq
From (\ref{Aelrelation}) and (\ref{A11pos}) it follows that $A_{21}$
and $A_{12}$ are nonzero and have opposite signs. If $A_{21}$ is positive, $A_{12}$ and
hence $\frac{A_{12}}{A_{11}}$ are negative, so entry $M^{(2)}_{12}$ is negative
and $M^{(2)}_{22}$ is positive, which is just the situation in the
numerical example shifted cyclically
one column to the right and with rows 1 and 2 interchanged.
Hence without loss of generality we can take
\beqa
  A_{21}  & < & 0, \label{A21neg} \\
  A_{12}  & > & 0. \label{A12pos}
\eeqa
It follows that the first row has the only positive entry in column 1
of $M^{(1)}$ and both constraint row entries in column 2 of
$M^{(2)}$ are positive.  Hence row 1 is the unique pivot candidate in
iteration 1.
%\cg{tie}{
There are two possible choices of pivot in column 2 of iteration 2.
We shall use the {\em largest pivot} rule
to break a tie\cg{. This rule chooses from the possible pivots the one
of largest magnitude, and is the best choice from the point of view
of numerical stability. }
{ which breaks a tie by
choosing the largest possible pivot, which would normally be the
choice from the point of view of numerical stability.} To simplify the
presentation we assume that if a tie remains after
applying this rule, then the pivot in row 1 is chosen. This second
tie-break rule therefore breaks the 2/6-cycle pattern if the pivot size
criterion does not determine the pivot row.
It follows that row 2
is the pivot choice in column 2 of iteration 2 if and only if
\beq
  \frac{1}{A_{11}}> \frac{A_{12}}{A_{11}} \iff A_{12} < 1. \label{A12lt1}
\eeq

We have therefore proved the following proposition.
\begin{prop} \label{PropRowSelect}
If the conditions of Proposition \ref{Prop26pattern} are met and row
selection ties are resolved by choosing the largest pivot and the
columns are selected in the 2/6-cycle order, then the necessary and
sufficient conditions for row 1 to be selected in odd iterations and
row 2 in even iterations are $0<A_{11}$ and $0<A_{12}<1$.
%, in
%which case $A_{21}<0$, or $0<A_{21}<1$, in which case $A_{12}<0$.
%Both condition yield the same examples delayed by one iteration.
%(Without loss of generality we assume that $A_{21}<0$ and $A_{12}>0$).
\end{prop}
The conditions guaranteeing that column 1 is chosen in
$M^{(1)}$ by the most negative reduced cost rule
rather than column 2 or 3 are
\beqa
  -1 & < & \mu, \label{mugtneg1} \\
  -1 & < & -(A_{11}+1)+\mu A_{21}
           \iff \mu < \frac{A_{11}}{A_{21}}. \label{T1c1vc3}
\eeqa
It follows from (\ref{A11pos}) and (\ref{A21neg}) that $\mu$ is negative.
Column 1 is guaranteed to be chosen rather than column 4 if and only if
\beqs
 -1 < -A_{12}-\mu A_{11} \iff \mu < \frac{1-A_{12}}{A_{11}},
\eeqs
which is always true as this bound is positive by (\ref{A11pos}) and
(\ref{A12lt1}).

In $M^{(2)}$, column 5 has a positive cost entry so is not a candidate. The
necessary and sufficient conditions for column 2 to be a candidate and
be guaranteed to be chosen rather that columns 3 or 4 are
\beqa
    \mu & < & -\frac{A_{12}}{A_{11}}, \label{T2c2cand} \\
   \mu & < & -\frac{(2+A_{11}+\frac{1}{A_{11}}+ \frac{A_{12}}{A_{11}})}
                      {1-A_{21}},  \label{T2c2vc3} \\
   \mu & < & -\frac{A_{12}(1+\frac{2}{A_{11}})}{A_{11}+1}.
      \label{T2c2vc4}
\eeqa
Comparing (\ref{T2c2cand}) and (\ref{T2c2vc4}) we see that (\ref{T2c2cand})
is redundant if
\beqa
    -\frac{A_{12}(1+\frac{2}{A_{11}})}{A_{11}+1}
     & < & -\frac{A_{12}}{A_{11}} \nonumber \\
\iff -A_{12}A_{11}-2A_{12} & < & -A_{11}A_{12}-A_{12} 
\iff -A_{12} < 0, \nonumber
\eeqa
which is true by (\ref{A12pos}).
Comparing (\ref{T2c2vc3}) and (\ref{T2c2vc4}), then using (\ref{A21neg}) and
then (\ref{Aelrelation}), we see that (\ref{T2c2vc3})
is redundant if
\beqa
    & -\frac{A_{12}(1+\frac{2}{A_{11}})}{A_{11}+1}
      <  -\frac{(2+A_{11}+\frac{1}{A_{11}}+ \frac{A_{12}}{A_{11}})}
                      {1-A_{21}} \nonumber \\
\iff & A_{12}(A_{11}+2)(1-A_{21})
      > (A_{11}+1)(2A_{11}+A^2_{11}+1+A_{12}) \nonumber \\
\iff & A_{12}(A_{11}+2-A_{21}(A_{11}+2)-A_{11}-1)
      > (A_{11}+1)^3 \nonumber \\
\iff & -A_{12}A_{21}(A_{11}+2)+A_{12}
      > (A_{11}+1)^3 \nonumber \\
\iff & (1+A_{11}+A^2_{11})(A_{11}+2)+A_{12}
      > (A_{11}+1)^3 \nonumber \\
\iff & (A_{11}+1)^3+1+A_{12}
      > (A_{11}+1)^3 \iff 1+A_{12} > 0, \nonumber
\eeqa
which (\ref{A12pos}) shows is true.
Comparing (\ref{T1c1vc3}) and (\ref{T2c2cand}), then using (\ref{A21neg})
and then
(\ref{Aelrelation}), we see that (\ref{T1c1vc3})
is redundant if
\beq
    -\frac{A_{12}}{A_{11}} < \frac{A_{11}}{A_{21}}
\iff  -A_{12}A_{21} > A^2_{11}
\iff  A^2_{11}+A_{11}+1 > A^2_{11}, \nonumber
\eeq
which (\ref{A11pos}) shows is true.

We have now shown that (\ref{T1c1vc3}), (\ref{T2c2cand}) and (\ref{T2c2vc3})
are redundant, so
(\ref{T2c2vc4}) is always the tightest upper bound. From this
and (\ref{mugtneg1}) it follows that
$\mu$ must lie in the range
\beq
   -1 < \mu < -\frac{A_{12}(A_{11}+2)}{A_{11}(A_{11}+1)},
       \label{mulimits}
\eeq
and there is a positive gap between these bounds if and only if
\beq
-1 <  -\frac{A_{12}(A_{11}+2)}{A_{11}(A_{11}+1)} \;\;\;
\iff A_{12} < A_{11}\left(\frac{A_{11}+1}{A_{11}+2}\right). \label{A12A11reln}
\eeq

If the left hand inequality in (\ref{mulimits}) is reversed, then
column 2 will be chosen rather than column 1 in $M^{(1)}$, and if the
right hand inequality is reversed, then column 4 will be chosen
instead of column 2 in $M^{(2)}$.  In either case the $2/6$-cycle
pattern will be broken. If either inequality in (\ref{mulimits}) holds
as an equality, then the most negative reduced cost rule does not
uniquely determine the column to enter the basis. To simplify
presentation we assume that when this occurs a choice is made
which breaks the $2/6$-cycle pattern.

\cg{We have now proved the following proposition.}{We have now shown}
\begin{prop}
Assume that the the most negative
reduced cost column selection rule and the largest pivot row
degeneracy tie breaking rule are used.
Then a 4 variable 2 constraint degenerate LP problem will have the
2/6-cycle pattern and cycle indefinitely
if and only if the conditions of Propositions
\ref{Prop26pattern} and \ref{PropRowSelect}
hold and in addition (\ref{mulimits}) holds
(which implies (\ref{A12A11reln})).
\end{prop}
%2/6-cycle example to have the unique pivot sequence (1,1) and (2,2), and
%therefore for the example to cycle indefinitely, is that $0<A_{11}$,
%$0<A_{12}<1$ and that (\ref{mulimits}) holds (which implies relation
%(\ref{A12A11reln})).

\begin{figure}
\epsfysize=9cm
\epsfbox{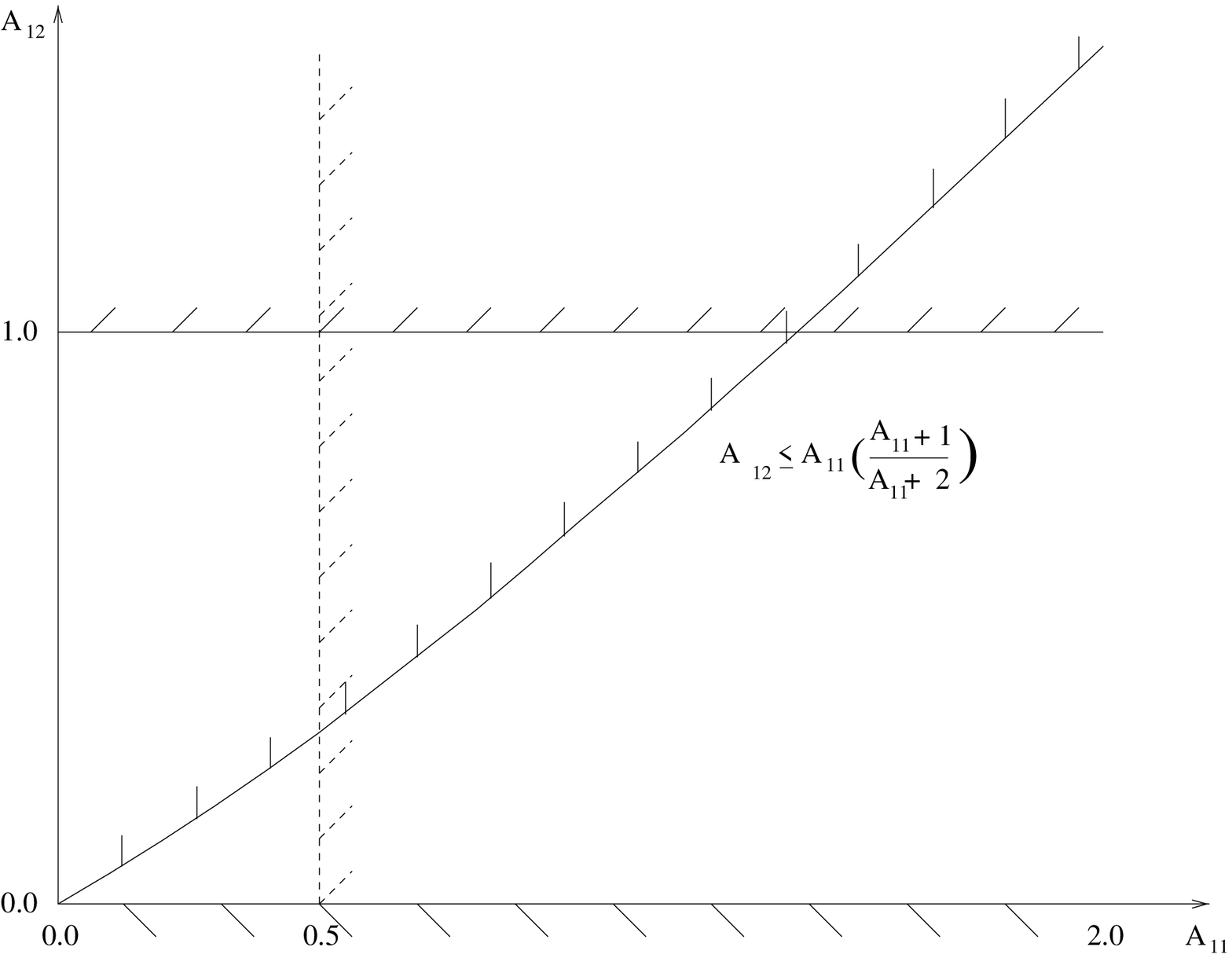}
\caption{Cycling region is unshaded. (Also cycles for \expand~if $A_{11}\le \half$)}
\label{cyclingregionfig}
\end{figure}

The unshaded area in Figure \ref{cyclingregionfig} (ignoring the dashed
constraint) shows the region where the problem cycles indefinitely.
Taking $A_{11}=0.4$, $A_{12}=0.2$ and $\mu=-2.15/2.3$ and then scaling
the objective row by $2.3$, produces the example given in Section
\ref{IntroEx}.

A similar analysis to that leading to Proposition \ref{Prop26pattern}
for the case of a $2/4$-cycle example shows that the cost row must be
zero, so such examples cannot cycle. It is also straightforward to
show that there can be no cycling examples with all pivots in the same
constraint row, so there can be no problems with a single constraint.
In the $2/6$-cycle examples $A_{12}$ and $A_{21}$ must have different
signs, so it follows from Table \ref{tab:CoefVals} that the even and
odd iterations cannot be the same.  \cg{Hence}{It therefore follows that}
the $2/6$-cycle examples are the simplest possible cycling examples.

\section{A cycling steepest-edge example} \label{SteepestEdge}

In the previous sections the column was selected using the original
Dantzig criterion of most negative reduced cost.  In the steepest-edge
method \cite{GoRe77} the column is selected on the basis of the most
negative ratio of the reduced cost to the length of the vector
corresponding to a unit change in the nonbasic variable. This normally
leads to a significant reduction in the number of iterations.
When steepest-edge column selection is used on the example in Section
\ref{IntroEx}, column 2 is chosen in $T^{(1)}$ instead of column 1 and in the
following iteration the problem is shown to be unbounded so the simplex
method terminates in 2 iterations.
%It can be shown that it is not
%possible to construct an example which cycles both for the most negative
%reduced cost and the steepest-edge criteria.
\cg{However by adding an extra row which affects the steepest-edge
weights but not the choice of pivot row, one can construct a steepest-edge
cycling example.}
{However one can construct a steepest-edge cycling example
by adding an extra row which affects the steepest-edge weights but not the choice
of pivot row. ??
However by adding an
extra row, which affects the steepest-edge weights but not the choice
of pivot row, a steepest-edge cycling example can be constructed.}

To preserve the 2/6-cycle pattern of the example, any extra
constraints must behave like the objective row in that they must
satisfy (\ref{baArel}). We shall now construct an example that has a
single candidate column in column 2 of $T^{(2)}$. We do this by
selecting $\mu$ so that the $x_4$ objective coefficient in $T^{(2)}$
is zero. It follows from Table \ref{tab:CoefVals} that the required value
is $\mu=-1.75$, and this results in the tableaux shown in
Table~\ref{steepestedgeTb}, omitting the third rows. Note that column
1 would not now be selected in $T^{(1)}$ either by the most negative
reduced cost criterion or by the steepest-edge criterion. We now
introduce a constraint that will leave the steepest-edge weight of
column 1 of $T^{(1)}$ unaltered but increase the weight of column
2. If the entries in this constraint are scaled up sufficiently, we
can make steepest-edge choose column 1. Using $a=[0,-20.0]$ and
applying (\ref{baArel}) we get the third row of tableau $T^{(1)}$.
We set the right-hand side of this constraint to 1,
\cg{which}{and this} ensures
that this constraint is not involved in any of the pivot
choices even when the matrix coefficients are perturbed by a small
amount. With this extra row added the steepest-edge reduced costs for
columns 1 and 2 of $T^{(1)}$ are $-0.127$ and $-0.087$, which leads to the
selection of column 1 as required.

\begin{table}
\caption{Cycling example with steepest-edge column selection}
\centerline{
\begin{tabular}{rrrrrrrrcrr}
$x_1$ & $x_2$ & $x_3$ & $x_4$ & $x_5$ & $x_6$ & $x_7$ & $I$ \\
\hline
   0.4 &   0.2 &  -1.4 &  -0.2 &  1.0  &     &     &     & = & 0 \\
  -7.8 &  -1.4 &   7.8 &   0.4 &       & 1.0 &     &     & = & 0 & ~~~~$T^{(1)}$ \\
   0.0 &  -20.0 & 156.0 &  8.0 &       &     & 1.0 &     & = & 1 \\
  -1.0 &  -1.75 &  12.25 &   0.5 &       &     &     & 1.0 & = & 0 \\
%  -0.8 &  -1.4 &   9.8 &   0.4 &       &     &     & 1.0 & = & 0 \\
\hline
   1.0 &   0.5  &  -3.5 &  -0.5 &   2.5 &     &     &    &= & 0 \\
       &   2.5  & -19.5 &  -3.5 &  19.5 & 1.0 &     &    &= & 0 & ~~~~~~~~$T^{(2)}$ \\
       &  -20.0 &  156.0 & 8.0 &   0.0 &      & 1.0 &    & = & 1 \\
       &   -1.25 &   8.75 &  0.0 &   2.5 &      &     & 1.0 & = & 0 \\
%       &   -1.0 &   7.0 &  0.0 &   2.0 &      &     & 1.0 & = & 0 \\
\hline
   1.0 &     &  0.4 &   0.2 &  -1.4 &  -0.2 &     &    & = & 0 \\
       & 1.0 & -7.8 &  -1.4 &   7.8 &   0.4 &     &    & = & 0 & ~~~~$T^{(3)}$ \\
       &     &   0.0 &  -20.0 & 156.0 &  8.0 & 1.0 &     & = & 1 \\
       &     & -1.0 &  -1.75 &  12.25 &   0.5  &     & 1.0 & = & 0 \\
%      &     &  -0.8 &  -1.4 &   9.8 &   0.4 &     & 1.0 & = & 0 \\
\hline
\end{tabular}
}
\label{steepestedgeTb}
\vspace{4mm}
\end{table}

\def\UserTol{u}

\section{Analysis of the \expand~procedure} \label{EXPAND}

The analysis given by Gill {\em et al} \cite{GMSW89} of their \expand~
procedure proves that the objective function can never return to a
value it had at a previous iteration. The \expand~procedure however
relaxes the constraints at each iteration, so the fact that the
objective function continually improves does not prove that the method
will not return to a previous basic solution. \cg{}{If this could occur the
method might cycle.} In Section \ref{EXPANDratiotest} we describe the
\expand~procedure and in Section \ref{EXPANDcycles} derive the
necessary and sufficient condition for cycling still to occur with the
2/6-cycle examples when using \expand. We do this by deriving an
expression for the values of every variable at every iteration, a
task that is made tractable by the special structure of the
$2/6$-cycle examples.

\subsection{The \expand~ratio test} \label{EXPANDratiotest}
The \expand~approach to resolving degeneracy is described by Gill~{\em et
al\/} in~\cite{GMSW89} for the general bounded LP problem.
The examples in this paper have single sided bounds and are of the form
\beqs
\begin{array}{lll}
{\rm minimize}\quad  & c^Tx \\
{\rm subject~to}\quad & Mx=b,&\quad x\ge 0.
\end{array}
\eeqs
For simplicity, \expand~is discussed here for this problem.
Assuming that all the variables are feasible ($x\ge 0$), the standard ratio
test for the simplex method determines the maximum step $\alpha$ in the
direction $p$ corresponding to the pivotal column such that the variables remain feasible,
that is $x-\alpha p\ge 0$.  For each $j$, the step which zeroes $x_j$ is
$\alpha_j=x_j/p_j$ if $p_j>0$, otherwise $\alpha_j=\infty$.  The maximum
feasible step is therefore $\alpha=\alpha_r=\min_j\alpha_j$ and the variable
to leave the basis is $x_r$.

\expand~is based on the use of an increasing primal feasibility tolerance
$\delta$.  During a particular `current' simplex iteration, this tolerance
has the value $\delta=\tilde\delta+\tau$, where $\tilde\delta$ was the value
of $\delta$ in the previous iteration. At the beginning of the current
iteration each variable satisfies its expanded bound $x_j\ge-\tilde\delta$.
Since $-\delta<-\tilde\delta$, it is always possible to ensure that
$\alpha>0$, so there is a strict decrease in the objective function.

The \expand~ratio test makes two passes through the entries in the pivotal
column $p$.
\begin{itemize}
\item The first pass determines the maximum acceptable step
$\alpha^{\rm max}>0$ so that each basic variable satisfies its new expanded bound $x_j\ge-\delta$.

\item
The second pass determines a variable $x_r$ to leave the basis. $x_r$
is the variable with the largest acceptable pivot and is defined by
\beqs
r=\arg \max_j p_j{\rm~such~that~} \alpha_j\le\alpha^{\rm max}{\rm~where~}
\alpha_j=\cases{x_j/p_j& $p_j>0$\cr \alpha_j=\infty&otherwise.}
\eeqs
Define $\alpha^{\rm full}=\alpha_r$. This is the step
necessary to zero $x_r$.  Note that if $x_r<0$
and $p_r>0$ then $\alpha^{\rm full}$ will be negative.

\item A minimum acceptable step
\beqs
\alpha^{\rm min}={\tau\over p_r}
\eeqs
is calculated. If $x_r=-\tilde\delta$ then this is the maximum step that
can be taken whilst maintaining feasibility with respect to the new expanded
bounds.

\item The actual step returned by the \expand~ratio test is
\beqs
\alpha=\max(\alpha^{\rm min},~\alpha^{\rm full}).
\eeqs
We refer to these two alternative step
sizes as the {\em min} and the {\em full} step.
\end{itemize}

The initial values of the nonbasic variables are zero. In the
2/6-cycle examples the initial values of the basic variables are also
zero. The initial value of the expanding feasibility tolerance is
denoted by $\tau u$, where $u\geq 0$, and the tolerance during
iteration $n$ is denoted by $\tau u^n$. It follows that
$u^n=u+n$.

\subsection{Conditions under which cycling occurs with the \expand~ratio test}
\label{EXPANDcycles}

In this section we analyse the behaviour of the $2/6$-cycle problems when using
the \expand~ratio test and derive necessary and sufficient conditions for
the $2/6$-cycle problems to cycle indefinitely.

%If the tolerance to which the problem is to be solved is expressed as
%$\delta_{\rm f}=2\UserTol\tau$ then the value of the expanding
%feasibility tolerance for the ratio test in iteration $n$ is
%$\tilde\delta=(\UserTol+n)\tau$.

The action of the \expand~ratio test depends on whether the iteration
number is even or odd, so we consider separately the behaviour in
iterations $n=2k+1$ and $n=2k+2$ for $k\ge0$. We assume that the
pivot columns are selected in the 2/6-cycle order and
derive necessary and sufficient conditions for \expand~to
select a pivot
in the first row in odd iterations and have a unique pivot
in the second row in even
iterations. We also show that the min step is taken when the
pivot is in row 1 and the full step is taken when the pivot is in row
2.

Let $x^n_j$ denote the value of $x_j$ at the start of iteration
$n$. The subscripts of $x$ are calculated modulo 6.

For iteration $2k+1$ the pivotal column is $\left[\matrix{A_{11}&
A_{21}}\right]^T$ and the values of the basic variables at the start
of the iteration are respectively $x_{2k-1}^{2k+1}$ and
$x_{2k}^{2k+1}$. Since $A_{21}<0$ and $A_{11}>0$, only $x_{2k-1}$ moves
towards its bound, so it is the sole candidate to leave the basis.  The
second pass of the \expand~ratio test returns
\beqs
\alpha^{\rm full}={x_{2k-1}^{2k+1}\over A_{11}},
\eeqs
and if
\beq
x_{2k-1}^{2k+1}\le\tau, \label{eq:Show1}
\eeq
the min step will be taken so
\beqs
\alpha = {\tau\over A_{11}}.
\eeqs
It follows that if (\ref{eq:Show1}) holds, the changes in variable
values are as given in row 1 of Table \ref{tab:AnyTwoIterations}.

For iteration $2k+2$ the pivotal column is
$\left[\matrix{A_{12}/A_{11}& 1/A_{11}}\right]^T$ and the values of
the basic variables at the start of the iteration are respectively
$x_{2k+1}^{2k+2}$ and $x_{2k}^{2k+2}$. Since $A_{11}>0$ and
$A_{12}>0$, both variables move towards their bound.  The first pass
of the \expand~ratio test returns
\beqs
\alpha^{\rm max} =\min\left({x_{2k+1}^{2k+2}+\tau u^{2k+2}\over A_{12}/A_{11}},
                {x_{2k}^{2k+2}+\tau u^{2k+2}\over 1/A_{11}} \right).
\eeqs
A sufficient condition for the pivot to be in row 2 is that $A_{12}<1$ and that
the pivot is acceptable. It is acceptable if and only if
\beqs
\begin{array}{lrcl}
&{\displaystyle {x_{2k}^{2k+2}\over 1/A_{11}}}&\le& \alpha^{\rm max} \\
\iff &A_{11}x_{2k}^{2k+2}&\le&{\displaystyle A_{11}\min\left({x_{2k+1}^{2k+2}+
\tau u^{2k+2}\over A_{12}}, x_{2k}^{2k+2}+\tau u^{2k+2} \right)}.
\end{array}
\eeqs
Clearly $x_{2k}^{2k+2}<x_{2k}^{2k+1}+\tau u^{2k+2}$, so the pivot in row
2 is acceptable if and only if
\beq
A_{12}x_{2k}^{2k+2} \leq x_{2k+1}^{2k+2} + \tau u^{2k+2}. \label{eq:Show2}
\eeq
Also, provided that
\beq
x_{2k}^{2k+2}\ge\tau\, \label{eq:Show3}
\eeq
then $\alpha^{\rm full}=A_{11}x_{2k}^{2k+2}\ge\alpha^{\rm min}=A_{11}\tau$,
so the full step
$\alpha^{\rm full}$ is taken and the \expand~ratio test returns
\beqs
\alpha=A_{11}x_{2k}^{2k+2}.
\eeqs
Hence if (\ref{eq:Show2}) and (\ref{eq:Show3}) hold, then the changes
in values are as given in row 2 of Table \ref{tab:AnyTwoIterations}.

\begin{sidewaystable}
%\begin{table}
\caption{Changes in values of variables over two iterations
~~~~~~~~~~~~~~~~~~~~~~~~~~~~~~~~~~~~~~~~~~~~~~~~~~~~~~~~~~~~~~~~~~~~~~~~~~~~~~~~~~~~~~~~~~~~~~~~}
\vspace{5mm}
%\centerline{
\begin{tabular}{c|lll|c}
$n$&\multicolumn1{c}{Entering}&\multicolumn1{c}{Leaving}&
 \multicolumn1{c}{Remaining}&\multicolumn1{c}{Step}\\ \hline
\vbox to 18pt {}
$2k+1$&
$x_{2k+1}^{2k+2}=x_{2k+1}^{2k+1}+{\displaystyle{\tau\over A_{11}}}$&
$x_{2k-1}^{2k+2}=x_{2k-1}^{2k+1}-\tau$&
$x_{2k  }^{2k+2}=x_{2k}^{2k+1}-\tau {\displaystyle{A_{21}\over A_{11}}}$&
Pivot row 1. Min step\\
\vbox to 14pt {}
$2k+2$&
$x_{2k+2}^{2k+3}=x_{2k}^{2k+2} A_{11}$&
$x_{2k  }^{2k+3}=0$&
$x_{2k+1}^{2k+3}=x_{2k+1}^{2k+2}-x_{2k}^{2k+2} A_{12}$&
Pivot row 2. Full step\\
\end{tabular}
%}
\label{tab:AnyTwoIterations}
\vspace{-5mm}
%\end{table}

{\vskip2cm}
\def\ds#1{\mbox{$\displaystyle #1$}}
%\begin{table}[h]
\caption{Expressions for the values of each variable over any two iterations.
 ~~~\ds{s_k = \sum_{i=0}^k A_{11}^i},~~\ds{S_k = \sum_{i=0}^k (k+1-i)A_{11}^i}.~~~~~~~~~~~~~~~~~~~ }
\vspace{5mm}
{\scriptsize
%\centerline{
\begin{tabular}{c|cccccc|c|c|}
%&Entering&{\multicolumn3{c}{Nonbasic}&Leaving&Basic\\
$n$&$x_{2k+1}^n$&$x_{2k+2}^n$&$x_{2k+3}^n$&$x_{2k+4}^n$&$x_{2k+5}^n$
&$x_{2k+6}^n$
 & Expanded & Normal \\\hline
\vbox to 14pt {}
& $-\tau S_{k-2}$ & $0$ & $-\tau S_{k-1}$ & $0$ & $\tau(1-S_k)$
   & $-\tau A_{21}s_{k-1}$ & & \\ \cline{2-9}
\vbox to 19pt {}
$2k+1$ & \ds{A_{11}} & \ds{A_{12}} & \ds{-(A_{11}+1)} & \ds{-A_{12}}
 & \ds{1} & \ds{0} 
 & \ds{\tau(1-S_k+u_{2k+1}){1\over A_{11}}} & \ds{\tau(1-S_k){1\over A_{11}}}\\
       & \ds{A_{21}} & \ds{-(A_{11}+1)} & \ds{-A_{21}} & \ds{A_{11}} 
 & \ds{0} & \ds{1}
 & \ds{\infty} & \ds{\infty} \\
\vbox to 1pt {} & \ds{\uparrow} & & & & & & &\\ \hline
\vbox to 19pt {}
& \ds{\tau({1\over A_{11}}- S_{k-2})} & \ds{0} & \ds{-\tau S_{k-1}} & \ds{0}
 & \ds{-\tau S_k} & \ds{-\tau{A_{21}\over A_{11}}s_{k}}
 & & \\  \cline{2-9}
\vbox to 19pt {}
\ds{2k+2} & \ds{1} & \ds{ {A_{12}\over A_{11}}} & \ds{-(1+{1\over A_{11}})}
 & \ds{-{A_{12}\over A_{11}}} & \ds{{1\over A_{11}}} & \ds{0}
 & \ds{\tau({1\over A_{11}}-S_{k-2}+u_{2k+2}){A_{11}\over A_{12}}}
 & \ds{\tau({1\over A_{11}}-S_{k-2}){A_{11}\over A_{12}}} \\
\vbox to 19pt {}
          & \ds{0} & \ds{{1\over A_{11}}} & \ds{{A_{21}\over A_{11}}} &
  \ds{-(1+{1\over A_{11}})} & \ds{-{A_{21}\over A_{11}}} & \ds{1}
 & \ds{\tau(-{A_{21}\over A_{11}}s_k+u_{2k+2})A_{11}}
 & \ds{-\tau A_{21}s_k}
\\
\vbox to 1pt {} & & \ds{\uparrow} & & & & & &\\ \hline
\vbox to 14pt {}
& \ds{\tau(1-S_{k+1})}  & \ds{-\tau A_{21}s_k} & \ds{-\tau S_{k-1}} & \ds{0} 
 & \ds{-\tau S_k} & \ds{0} & &\\ \hline
\end{tabular}}
%}
\label{tab:FirstSevenIterations}
%\end{table}
\end{sidewaystable}

From the changes in the values of variables given in Table
\ref{tab:AnyTwoIterations}, the expressions in Table
\ref{tab:FirstSevenIterations} for the values of each
variable over any two iterations are established by induction.
To simplify notation we introduce the quantities
%In Table \ref{tab:FirstSevenIterations} and below the notation
$s_k$ and $S_k$ defined by
\beqas
s_k  =  \sum_{i=0}^k A_{11}^i, & \;\;\;
S_k  =  \sum_{i=0}^k(k+1-i)A_{11}^i, &\mathrm{~for~all~} k \geq 0, \\
s_k  = 0, &\;\;\; S_k = 0, &\mathrm{~for~all~} k<0.
\eeqas
Note that since $A_{11}>0$,
$s_k$ and $S_k$ are nonnegative. Also
\beqa
S_k - S_{k-1}  = s_k, \mathrm{~for~ all~} k, \nonumber \\
s_k = 1+A_{11}s_{k-1}, \mathrm{~for~all~} k \geq 0. \label{srelation}
\eeqa
The expressions in Table \ref{tab:FirstSevenIterations} allow condition
(\ref{eq:Show2}) to be expressed as $G_k\ge0$, where
$G_k$ for $k \geq 0$ is defined by
\beqas
G_k &=& {A_{12}A_{21}\over A_{11}}s_k + {1\over A_{11}} - S_{k-2} + u^{2k+2}.
\eeqas
A necessary and sufficient condition on $A_{11}$ for $G_k\ge0$ is
established by considering
\beqas
\Delta G_k &\equiv& G_{k+1}-G_{k}  \\
           &=& {A_{12}A_{21}\over A_{11}}(s_{k+1}-s_{k})
               -(S_{k-1}-S_{k-2}) + u^{2k+4} - u^{2k+2}  \\
           &=& -{1+A_{11}+A_{11}^2\over A_{11}}A_{11}^{k+1}
                - s_{k-1} + 2 \\
           &=& -(s_{k-1}+A_{11}^k+A_{11}^{k+1}+A_{11}^{k+2})+2 \\
           &=& -s_{k+2} + 2.
\eeqas
It follows that $\Delta G_k \geq 0 \iff s_{k+2} \leq 2$.  If
$0<A_{11}\leq {1\over 2}$, then $s_{k+2}$ increases to a limit
$s_{\infty}$, where $s_{\infty} \leq 2$. In this case $\Delta G_k \geq
0$ for all $k$ so $G_{k+1} \geq G_{k}$ for all $k$, and also $G_0 \geq
u+{1\over 2} > 0$, so $G_k>0$ for all $k\geq 0$. If $A_{11} > {1\over
2}$, then there exists an $\epsilon$ and $K$ such that
$s_{k+2}>2+\epsilon$ for all $k\geq K$. It follows that for suitably
large $k$, $G_k<0$. Hence for positive $A_{11}$ the necessary and
sufficient conditions for $G_k$ to be nonnegative for all $k$ is that
$A_{11} \leq {1\over 2}$.

\begin{prop}
Assume that the conditions of Proposition \ref{Prop26pattern} are met
and the \expand~row selection method is used and the columns are
selected in the 2/6-cycle order. Then necessary and sufficient conditions
for cycling to occur are that $0< A_{11}\leq {1\over 2}$ and $0<A_{12}<1$.
\label{A11leHalfProp}
\end{prop}

\noindent
{\bf Proof:}

%\vspace{2mm}
%\noindent
{\bf Sufficient conditions:}
%\vspace{2mm}

%We first prove the sufficient conditions.
We show by induction that the values of the variables at the
start of odd iterations are as given in Table
\ref{tab:FirstSevenIterations} and that these values
lead to the correct choice of pivot row for the 2/6-cycle pattern.

Initially all the variables have the value zero, so $x_j^1=0$. Hence
the values in Table \ref{tab:FirstSevenIterations} are correct for
$n=1$.  Assume now that for some $k$ the values in Table
\ref{tab:FirstSevenIterations} are correct at the start of iteration
$2k+1$.

In iteration $2k+1$, since $s_k$ is non-negative,
$x_{2k-1}^{2k}\le\tau$, so (\ref{eq:Show1}) holds and it follows that the
changes in the values of variables are as given by row 1 of Table
1. From this and (\ref{srelation}) we get
\beqa
x_{2k+6}^{2k+2} &=&-\tau A_{21}s_{k-1} - \tau{A_{21}\over A_{11}} \\ \nonumber
  &=&-\tau{A_{21}\over A_{11}}(A_{11}s_{k-1} + 1) \\ \nonumber
  &=& -\tau{A_{21}\over A_{11}}s_k. \nonumber
\eeqa
All the other values are straightforward, so we have deduced the values
given in Table \ref{tab:FirstSevenIterations} at the start of iteration
$2k+2$.

Substituting these values into (\ref{eq:Show2}) we see that the pivot in row 2 is
acceptable if and only if
\beqs
 -{A_{12}A_{21}\over A_{11}}s_k \leq {1\over A_{11}} - S_{k-2} + u^{2k+2},
\eeqs
which is true provided $A_{11}\leq {1\over 2}$.
%\beqs
%\begin{array}{rcl}
%\Delta_k&=&A_{12}x_{2k}^{2k+1}-x_{2k+1}^{2k+1}-(\UserTol+2k+2)\tau
%\\\noalign{\medskip\noindent so\medskip}
%\Delta_k&=&\Delta_{k-1}+\delta_k,
%\\\noalign{\medskip\noindent where\medskip}
%\delta_k
%&=&\Delta_k-\Delta_{k-1}\\ \noalign{\smallskip}
%&=&A_{12}(x_{2k}^{2k+1}-x_{2k-2}^{2k-1})-(x_{2k+1}^{2k+1}-x_{2k-1}^{2k-1})-2\tau  \\ \noalign{\smallskip}
%&=&\displaystyle{-\tau A_{12}A_{21}A_{11}^{k-1}+\tau\sum_{i=0}^{k-2}A_{11}^i-2\tau} \\ \noalign{\smallskip}
%&=&\displaystyle{\tau (1+A_{11}+A_{11}^2)A_{11}^{k-1}+\tau\sum_{i=0}^{k-2}A_{11}^i-2\tau}  \\ \noalign{\smallskip}
%&=&\displaystyle{\left(\sum_{i=0}^{k+1}A_{11}^i-2\right)\tau}.
%\end{array}
%\eeqs
%Since~(\ref{eq:Show2}) is assumed to hold for $k-1$ it follows that $\Delta_{k-1}\le0$ so~(\ref{eq:Show2}) holds for $k$ if $\delta_k\le0$. A sufficient
%condition for $\delta_k\le0$ independent of $k$ is that $0<A_{11}\le\frac{1}{2}$.
Also
\beqs
x_{2k}^{2k+2}=-\tau {A_{21}\over A_{11}}\sum_{i=0}^{k} A_{11}^i\geq -\tau {A_{21}\over A_{11}}= -\tau{1+A_{11}+A_{11}^2\over A_{12}}>\tau,
\eeqs
since $A_{11}>0$~(\ref{A11pos}) and $A_{12}<1$. Hence~(\ref{eq:Show3}) holds, so the changes in the values of variables are as given in row
2 of Table \ref{tab:AnyTwoIterations}. The new value for $x_{2k+1}$ is given by
\beqa
x_{2k+1}^{2k+3} &=& \tau\left({1\over A_{11}} - S_{k-2}
                    + {A_{12}A_{21}\over A_{11}}s_k\right) \nonumber \\
   &=& \tau\left({1\over A_{11}} - S_{k-2}
                    -{1\over A_{11}}s_k-s_k-A_{11}s_k\right) \nonumber \\
   &=& \tau\left({1\over A_{11}} - S_{k-2}
                    -({1\over A_{11}}+s_{k-1})-s_k-(s_{k+1}-1)\right) \nonumber \\
   &=& \tau(1-S_{k+1}), \nonumber
\eeqa
%\beqs
%\begin{array}{rcl}
%x_{2k+1}^{2k+3}
%&=&{\displaystyle{\tau\over A_{11}} -\tau\sum_{i=0}^{k-2}(k-1-i)A_{11}^i+\tau {A_{12}A_{21}\over A_{11}}\sum_{i=0}^{k}A_{11}^i}\\\noalign{\smallskip}
%&=&{\displaystyle{\tau\over A_{11}} -\tau\sum_{i=0}^{k-2}(k-1-i)A_{11}^i-\tau(1+A_{11}+A_{11}^2)\sum_{i=0}^{k}A_{11}^{i-1}}\\\noalign{\smallskip}
%&=&{\displaystyle\tau               -\tau\sum_{i=0}^{k-2}(k-1-i)A_{11}^i-3\tau\sum_{i=0}^{k}A_{11}^i-2\tau A_{11}^{k+1}
%                     -\tau A_{11}^{k+2}}\\\noalign{\smallskip}
%&=&{\displaystyle\tau               -\tau\sum_{i=0}^{k+1}(k+2-i)A_{11}^i}, \\\noalign{\smallskip}
%\end{array}
%\eeqs
which is the value given in Table~\ref{tab:FirstSevenIterations}.
All the other values at the start
of iteration $2k+3$ follow straightforwardly and are as shown in
Table~\ref{tab:FirstSevenIterations}. These values are the values in
Table~\ref{tab:FirstSevenIterations} for $k$, with the $k$
replaced by $k+1$. This completes the induction and shows that the 2/6-cycle
pattern continues indefinitely.

\vspace{2mm}
%\noindent
{\bf Necessary conditions:}
%\vspace{2mm}

As discussed in Section \ref{2/6form}, $A_{11}>0$ and we can choose
$A_{12}>0$, in which case $A_{21}<0$. Since $x_5^1=0$, the first
iteration takes the min step and so $x_1^2=\tau/A_{11}$. The pivot in
row 1 in iteration 2 is acceptable if
\beqa
&{\tau\over A_{11}}{A_{11}\over A_{12}} & \leq \alpha^{\rm max} \nonumber \\
   \iff & {\tau\over A_{12}}
       & \leq -{\tau A_{21}\over A_{11}} + \tau(u+2) \nonumber \\
   \iff & 1 & \leq {1\over A_{11}}+1+A_{11}+u+2,
\eeqa
which is true. Hence if $A_{12}>1$, the pivot will be in row 1 in
iteration 2 and the 2/6-cycle pattern will be broken.  If
$A_{11}>{1\over 2}$, then the argument prior to
Proposition~\ref{A11leHalfProp} shows there is
a first value of $k$, $\hat{K}$ say, such
that $G_{\hat{K}}<0$. As shown above, for all $k < \hat{K}$ the 2/6-cycle pattern is
maintained and the variable values are as in
Table~\ref{tab:FirstSevenIterations}. Therefore in iteration
$2\hat{K}$ the pivot in row 2 is not acceptable, so the pivot must
be in row 1. This breaks the 2/6-cycle pattern.  $\; \Box$
\vspace{3mm}

The conditions derived in Section \ref{2/6form} for the minimum reduced cost
criterion to choose pivot columns in the 2/6-cycle pattern relied on
the conditions $A_{11}>0$ and $0<A_{12}<1$. These
conditions have been established in Proposition \ref{A11leHalfProp}
for the case of
\expand~row selection, so it follows that (\ref{mulimits})
and (\ref{A12A11reln}) still hold. From (\ref{A12A11reln}) and the
fact that $A_{11}\leq {1\over 2}$ it follows that $A_{12}<{3\over
10}$, which is tighter that $A_{12}<1$, which is therefore redundant.
We have now shown the following proposition.

\begin{prop}
A 4 variable 2 constraint degenerate LP problem will have the
2/6-cycle pattern and cycle indefinitely when using the most negative
reduced cost column selection rule and the \expand~row selection
rule if and only if the conditions of Proposition
\ref{Prop26pattern} hold and in addition $0<A_{11}\le\half$,
$0<A_{12}$ and (\ref{mulimits}) holds
(which implies relation (\ref{A12A11reln})).
\end{prop}

The shaded area in Figure \ref{cyclingregionfig} including the
$A_{11}\le\half$ constraint is the region where cycling occurs
when using \expand.  Note that the constraint $A_{12}<1$ is
now redundant.  Also note that in the example (\ref{prob1}), $A_{11} = 0.4 <
0.5$, so that example also cycles with \expand.

Finally note that the only way that \expand~can escape from the
$2/6$-cycle pattern is for it to select the first row as pivot row
in an even
iteration, and if this occurs the resulting tableau has the form

%\vspace{4mm}
\centerline{
\begin{tabular}{cccccc}
$x_1$ & $x_2$ & $x_3$ & $x_4$ & $x_5$ & $x_6$ \\ \noalign{\smallskip}
\hline
 $\frac{A_{12}}{A_{11}}$  & 1 & $-\frac{A_{12}}{A_{11}}-\frac{1}{A_{12}}$ & $-1$
  &   $\frac{1}{A_{12}}$ & $0$  \\  \noalign{\smallskip}
 $-\frac{1}{A_{12}}$  & $0$ &  $\frac{A_{21}}{A_{11}}+\frac{1}{A_{12}}+\frac{1}{A_{11}A_{12}}$ & $-1$ & $-\frac{A_{21}}{A_{11}}-\frac{1}{A_{11}A_{12}}$
  &   1 \\  \noalign{\smallskip}
 $-1-\mu{\Aaa\over \Aab}$  & 0 &
$*$ & $*$ & $-\frac{\mu}{\Aab}$ & 0\\  \noalign{\smallskip}
\hline
\end{tabular}
}
The constraint entries in the third and fourth columns are all
negative and the objective function coefficients in all except these
columns are nonnegative.  Since the problem is unbounded we
cannot be at an optimum, so one of these columns must be chosen.
The next iteration will then produce an unbounded step and
the method will terminate.
%
%It is easy to show that $0<A_{11}\le\frac{1}{2}$ is
%also a necessary condition for~\expand~to cycle.  Suppose that
%$A_{11}>\frac{1}{2}$. Then there exists an integer $K$ such that
%$\delta_k\ge\delta>0$ for all $k\ge K$. If $\Delta_K=\Delta$ then
%$\Delta_k\ge\Delta+(k-K)\delta$ for all $k\ge K$, in which case
%$\Delta_k>0$ for $k$ sufficiently large and (\ref{eq:Show2}) no longer
%holds. The pivoting sequence is no longer maintained so cycling
%terminates and degeneracy is resolved.

The above results are independent of the \expand~parameters $u$ and
$\tau$.  In \cite{GMSW89} it is suggested that the initial tolerance
$\tau u$ be taken to be half of the feasibility tolerance $\delta_{\rm
f}$ to which the problem is to be solved.  The value of $\tau$ is
chosen so that after a large number of iterations (typically
$K=10000$) the expanded tolerance approaches $\delta_{\rm f}$, at
which stage $\delta$ is reset to its original value $\delta_{\rm i}$. If
this is done with the 2/6-cycle examples after an even iteration, then
the problem returns to its initial state. If it is done after an odd
iteration then it returns to the even iteration case but with the
values all zero. It can be shown that in this case too the problem
cycles, so that in neither case does resetting break the cycle
pattern.

\section{Conclusions}

We have derived a three-parameter class of linear programming examples
which cause the simplex method to cycle indefinitely.  When written in
standard form, these examples have two constraints and 6 variables and
the coefficient pattern repeats every two iterations. These are the
simplest possible examples for which the simplex algorithm cycles.  We
have derived 4 inequalities between the parameters and shown that
these are the necessary and sufficient conditions for members of this
class to cycle with Dantzig's form of the simplex method.  We have shown how
to extend the examples so that they also cycle when the steepest-edge
column selection criterion is used. By adding the single bound,
$A_{11}\le \half$, we were able to characterise the examples that
also cycle using the
\expand~row-selection mechanism. This shows that despite the
fact that in the \expand~method the objective function is guaranteed to
improve each iteration,
the method is not guaranteed to prevent cycling.
The cycling behaviour is independent of the \expand~tolerance
parameters.  The bound $A_{11}\le \half$ is the only extra condition
that had to be applied to ensure that an example which would cycle
under the usual Dantzig rule, with largest pivot as the tie-breaker,
would also cycle using \expand. This extra bound does reduce somewhat
the number of cases that cycle, but does not eliminate the
problem.
\cg{Also the reduction is for problems where the degeneracy is exact.
For problems which are close to degenerate \expand~may cycle whereas
the original simplex method in exact arithmetic will not.}
{The reduction also has to be set against the fact that
\expand~can cycle when the degeneracy is not exact.}

All the coefficients in the examples (not just the 3
parameters) may be perturbed simultaneously by any small amount
without destroying the cycling behaviour.  The
2/6-cycle examples are therefore just points in a full dimensional set of
counter-examples, so there is a positive probability of encountering
cycling in randomly generated degenerate examples. In practice therefore
the \expand~procedure cannot be relied upon to prevent
cycling.
Provided we stay within the class of the
degenerate problems (\emph{i.e.}\ keep the right-hand side 0)
it is possible to vary the other coefficients by a significant
amount. Indeed we have constructed examples where the values are
totally different every 2 iterations and yet indefinite cycling
still occurs with \expand.

The examples have been tested on our own implementation of \expand~and
using {\sc minos 5.4}, which was written by the authors of \expand. In
both cases if no preprocessing is done the examples cycle
indefinitely.  {\sc minos} periodically does a reset operation (by
default after 10000 iterations). This returns the problem to its
initial state so cycling is still indefinite.

{\sc osl} \cite{ForTom95} uses some techniques from \expand.
%, however
%there is no evidence that it does the small perturbation steps each
%iteration.
In the examples in this paper, {\sc osl} 2.0 without
scaling or preprocessing and with Dantzig pricing cycles for 30 iterations
\cg{before reporting doing a perturbation,
but it then continues to cycle indefinitely.}
{before introducing a large perturbation, which resolves the degeneracy.}
{\sc cplex} 4.0.7
without scaling or preprocessing cycles for 400 iterations before
resolving the degeneracy by introducing a large perturbation.
{\sc xp}ress{\sc mp} 7.14 without scaling and with an even invert frequency
cycles indefinitely.
However the \cg{}{robust}
{\sc bqpd} code of Fletcher \cite{FL93}
detects degeneracy at the start of the first iteration, changes to the
dual, does one pivot, then finds that the
dual is infeasible.  This gives an improving direction in the primal,
which resolves the degeneracy. Finally, {\sc bqpd} detect unboundedness in
this direction and terminates, having done one pivot in total.

\bibliography{ref}
\bibliographystyle{abbrv}

%\input{appendix}

%\newpage

%Questions
%\begin{enumerate}

%\item Do we want to keep the 2/6 cycle name?

%\item Is it possible to do a dashed \verb+\+hline?

%\item Sort out $<$ or $\le$ questions.

%\item vr in row question.

%\item Should we put down the EPSRC grant so that we can cite this
%paper in the final report?

%\item What version of OSL do we have?

%\item shall we say that \expand works just as well on average with
%expansion = 0, or does that detract for the analysis in the paper.

%\item is the a better way of doing the $A^3=I$ analysis?

%\end{enumerate}

\end{document}